# Some Characterizations of Mannheim Partner Curves in Minkowski 3-space $E_1^3$

Tanju Kahraman[a], Mehmet Önder[a], Mustafa Kazaz[a], H. Hüseyin Uğurlu[b]

[a] Celal Bayar University, Department of Mathematics, Faculty of Arts and Sciences, , Manisa, Turkey;
tanju.kahraman@bayar.edu.tr, mehmet.onder@bayar.edu.tr, mustafa.kazaz@bayar.edu.tr,

[b] Gazi University, Gazi Faculty of Education, Department of Secondary Education Science and Mathematics Teaching, Mathematics Teaching Program, Ankara, Turkey; hugurlu@gazi.edu.tr

**Abstract:** In this paper, we give the characterizations of Mannheim Partner Curves in Minkowski 3-space $E_1^3$. Firstly, we classify these curves in $E_1^3$. Next, we give some relationships characterizing these curves and we show that Mannheim theorem is not valid for Mannheim partner curves in $E_1^3$. Moreover, by considering the spherical indicatrix of the Frénet vectors of those curves, we obtain some new relationships between the curvatures and torsions of the Mannheim partner curves in $E_1^3$.

## 1. INTRODUCTION

In the differential geometry, special curves have an important role. Especially, the partner curves, i.e., the curves which are related each other at the corresponding points, have drawn attention of many mathematicians so far. The well-known of the partner curves is Bertrand curves which are defined by the property that at the corresponding points of two space curves the principal normal vectors are common. Bertrand partner curves have been studied in ref. [1,2,3,4,13,15]. Ravani and Ku have transported the notion of Bertrand curves to the ruled surfaces and called Bertrand offsets[12]. Recently, Liu and Wang have defined a new curve pair for space curves. They called these new curves as Mannheim partner curves: Let $x$ and $x_1$ be two curves in the three dimensional Euclidean space $E^3$. If there exists a correspondence between the space curves $x$ and $x_1$ such that, at the corresponding points of the curves, the principal normal lines of $x$ coincides with the binormal lines of $x_1$, then $x$ is called a Mannheim curve, and $x_1$ is called a Mannheim partner curve of $x$. The pair $\{x, x_1\}$ is said to be a Mannheim pair. They showed that the curve $x_1(s_1)$ is the Mannheim partner curve of the curve $x(s)$ if and only if the curvature $\kappa_1$ and the torsion $\tau_1$ of $x_1(s_1)$ satisfy the following equation

$$\dot{\tau} = \frac{d\tau}{ds_1} = \frac{\kappa_1}{\lambda}(1 + \lambda^2 \tau_1^2)$$

for some non-zero constants $\lambda$. They also studied the Mannheim partner curves in the Minkowski 3-space and obtained the necessary and sufficient conditions for the Mannheim partner curves in $E_1^3$ [See 5 and 14 for details]. Moreover, Oztekin and Ergut[11] studied the null Mannheim curves in the same space. Orbay and Kasap gave[10] new characterizations of Mannheim partner curves in Euclidean 3-space. They also studied[9] the Mannheim offsets of ruled surfaces in Euclidean 3-space. The corresponding characterizations of Mannheim offsets of timelike and spacelike ruled surfaces have been given by Onder and et al [6,7].

In this paper, we give the new characterizations of Mannheim partner curves in Minkowski 3-space $E_1^3$. Furthermore, we show that the Mannheim theorem is not valid for Mannheim partner curves in $E_1^3$. Moreover, we give some new characterizations of the Mannheim partner curves by considering the spherical indicatrix of some Frénet vectors of the curves.

## 2. PRELIMINARIES

Minkowski 3-space $E_1^3$ is the real vector space $E^3$ provided with the standart flat metric given by

$$\langle,\rangle = -dx_1^2 + dx_2^2 + dx_3^2$$

where $(x_1, x_2, x_3)$ is a rectangular coordinate system of $E_1^3$. According to this metric, in $E_1^3$ an arbitrary vector $\vec{v} = (v_1, v_2, v_3)$ can have one of three Lorentzian causal characters; it can be spacelike if $\langle \vec{v}, \vec{v} \rangle > 0$ or $\vec{v} = 0$, timelike if $\langle \vec{v}, \vec{v} \rangle < 0$ and null(lightlike) if $\langle \vec{v}, \vec{v} \rangle = 0$ and $\vec{v} \neq 0$ [8]. Similarly, an arbitrary curve $\vec{\alpha} = \vec{\alpha}(s)$ can locally be spacelike, timelike or null (lightlike), if all of its velocity vectors $\vec{\alpha}'(s)$ are spacelike, timelike or null (lightlike), respectively. We say that a timelike vector is future pointing or past pointing if the first compound of the vector is positive or negative, respectively. For the vectors $\vec{x} = (x_1, x_2, x_3)$ and $\vec{y} = (y_1, y_2, y_3)$ in $E_1^3$, the vector product of $\vec{x}$ and $\vec{y}$ is defined by

$$\vec{x} \wedge \vec{y} = \begin{vmatrix} e_1 & -e_2 & -e_3 \\ x_1 & x_2 & x_3 \\ y_1 & y_2 & y_3 \end{vmatrix} = (x_2 y_3 - x_3 y_2, x_1 y_3 - x_3 y_1, x_2 y_1 - x_1 y_2)$$

where

$$\delta_{ij} = \begin{cases} 1 & i = j, \\ 0 & i \neq j, \end{cases} \quad e_i = (\delta_{i1}, \delta_{i2}, \delta_{i3}) \text{ and } e_1 \wedge e_2 = -e_3, \ e_2 \wedge e_3 = e_1, \ e_3 \wedge e_1 = -e_2.$$

The Lorentzian sphere and hyperbolic sphere of radius $r$ and center 0 in $E_1^3$ are given by

$$S_1^2 = \left\{ \vec{x} = (x_1, x_2, x_3) \in E_1^3 : \langle \vec{x}, \vec{x} \rangle = r^2 \right\}$$

and

$$H_0^2 = \left\{ \vec{x} = (x_1, x_2, x_3) \in E_1^3 : \langle \vec{x}, \vec{x} \rangle = -r^2 \right\},$$

respectively[6,7].

Denote by $\{\vec{T}, \vec{N}, \vec{B}\}$ the moving Frénet frame along the curve $\alpha(s)$ in the Minkowski space $E_1^3$. For an arbitrary spacelike curve $\alpha(s)$ in the space $E_1^3$, the following Frénet formulae are given,

$$\begin{bmatrix} \vec{T}' \\ \vec{N}' \\ \vec{B}' \end{bmatrix} = \begin{bmatrix} 0 & k_1 & 0 \\ -\varepsilon k_1 & 0 & k_2 \\ 0 & k_2 & 0 \end{bmatrix} \begin{bmatrix} \vec{T} \\ \vec{N} \\ \vec{B} \end{bmatrix} \quad (1_a)$$

where $g(\vec{T}, \vec{T}) = 1$, $g(\vec{N}, \vec{N}) = \varepsilon = \pm 1$, $g(\vec{B}, \vec{B}) = -\varepsilon$, $g(\vec{T}, \vec{N}) = g(\vec{T}, \vec{B}) = g(\vec{N}, \vec{B}) = 0$ and $k_1$ and $k_2$ are curvature and torsion of the spacelike curve $\alpha(s)$, respectively. Here, $\varepsilon$ determines the kind of spacelike curve $\alpha(s)$. If $\varepsilon = 1$, then $\alpha(s)$ is a spacelike curve with spacelike principal normal $\vec{N}$ and timelike binormal $\vec{B}$. If $\varepsilon = -1$, then $\alpha(s)$ is a spacelike curve with timelike principal normal $\vec{N}$ and spacelike binormal $\vec{B}$. Furthermore, for a timelike curve $\alpha(s)$ in the space $E_1^3$, the following Frénet formulae are given as follows,

$$\begin{bmatrix} \vec{T}' \\ \vec{N}' \\ \vec{B}' \end{bmatrix} = \begin{bmatrix} 0 & k_1 & 0 \\ k_1 & 0 & k_2 \\ 0 & -k_2 & 0 \end{bmatrix} \begin{bmatrix} \vec{T} \\ \vec{N} \\ \vec{B} \end{bmatrix}. \tag{$1_b$}$$

where $g(\vec{T},\vec{T}) = -1$, $g(\vec{N},\vec{N}) = g(\vec{B},\vec{B}) = 1$, $g(\vec{T},\vec{N}) = g(\vec{T},\vec{B}) = g(\vec{N},\vec{B}) = 0$ and $k_1$ and $k_2$ are curvature and torsion of the timelike curve $\alpha(s)$ respectively[See 8 and 16 for details].

**Definition 2.1. i) Hyperbolic angle:** *Let $\vec{x}$ and $\vec{y}$ be future pointing (or past pointing) timelike vectors in $E_1^3$. Then there is a unique real number $\theta \geq 0$ such that $<\vec{x},\vec{y}> = -|\vec{x}||\vec{y}|\cosh\theta$. This number is called the hyperbolic angle between the vectors $\vec{x}$ and $\vec{y}$.*

**ii) Central angle:** *Let $\vec{x}$ and $\vec{y}$ be spacelike vectors in $E_1^3$ that span a timelike vector subspace. Then there is a unique real number $\theta \geq 0$ such that $<\vec{x},\vec{y}> = |\vec{x}||\vec{y}|\cosh\theta$. This number is called the central angle between the vectors $\vec{x}$ and $\vec{y}$.*

**iii) Spacelike angle:** *Let $\vec{x}$ and $\vec{y}$ be spacelike vectors in $E_1^3$ that span a spacelike vector subspace. Then there is a unique real number $\theta \geq 0$ such that $<\vec{x},\vec{y}> = |\vec{x}||\vec{y}|\cos\theta$. This number is called the spacelike angle between the vectors $\vec{x}$ and $\vec{y}$.*

**iv) Lorentzian timelike angle:** *Let $\vec{x}$ be a spacelike vector and $\vec{y}$ be a timelike vector in $E_1^3$. Then there is a unique real number $\theta \geq 0$ such that $<\vec{x},\vec{y}> = |\vec{x}||\vec{y}|\sinh\theta$. This number is called the Lorentzian timelike angle between the vectors $\vec{x}$ and $\vec{y}$* [6,7].

In this paper, we study the Mannheim partner curves in $E_1^3$. We obtain the relationships between the curvatures and torsions of the Mannheim partner curves with respect to each other. Using these relationships, we give Mannheim's theorem for the Mannheim partner curves in Minkowski 3-space $E_1^3$.

## 3. MANNHEIM PARTNER CURVES IN MINKOWSKI 3-SPACE $E_1^3$

In this section, by considering the Frénet frames, we give the characterizations of Mannheim partner curves in Minkowski 3-space $E_1^3$.

**Definition 3.1.** *Let $C$ and $C^*$ be two curves in Minkowski 3-space $E_1^3$ given by the parametrizations $\alpha(s)$ and $\alpha^*(s^*)$, respectively, and let they have at least four continuous derivatives. If there exists a correspondence between the space curves $C$ and $C^*$ such that the principal normal lines of $C$ coincides with the binormal lines of $C^*$ at the corresponding points of curves, then $C$ is called as a Mannheim curve and $C^*$ is called as a Mannheim partner curve of $C$. The pair of $\{C, C^*\}$ is said to be a Mannheim pair*[5].

By considering the Lorentzian casual characters of the curves, from Definition 3.1, it is easily seen that there are five different types of the Mannheim partner curves in Minkowski 3-space $E_1^3$. Let the pair $\{C, C^*\}$ be a Mannheim pair. Then according to the characters of the curves $C$ and $C^*$ we have the followings:

**Case 1.** *The curve $C^*$ is timelike.*
If the curve $C^*$ is timelike, then there are two cases.

    **i)** The curve $C$ is a spacelike curve with a timelike principal normal. In this case, we say that the pair $\{C, C^*\}$ is a Mannheim pair of the type 1.

    **ii)** The curve $C$ is a timelike curve. In this case, we say that the pair $\{C, C^*\}$ is Mannheim pair of the type 2.

**Case 2.** *The curve $C^*$ is spacelike.*
If the curve $C^*$ is a spacelike curve, then there are three cases;

    **iii)** The curve $C^*$ is a spacelike curve with a timelike binormal vector and the curve $C$ is a spacelike curve with a timelike principal normal vector. In this case, we say that the pair $\{C, C^*\}$ is a Mannheim pair of the type 3.

    **iv)** The curve $C^*$ is a spacelike curve with a timelike binormal vector and the curve $C$ is a timelike curve. In this case, we say that the pair $\{C, C^*\}$ is a Mannheim pair of the type 4.

    **v)** The curve $C^*$ is a spacelike curve with a timelike principal normal vector and the curve $C$ is a spacelike curve with a timelike binormal vector. In this case, we say that the pair $\{C, C^*\}$ is a Mannheim pair of the type 5.

**Theorem 3.1.** *The distance between corresponding points of the Mannheim partner curves is constant in $E_1^3$.*

*Proof.* Let consider the case that the pair $\{C, C^*\}$ is a Mannheim pair of the type 1. From Definition 3.1, we can write
$$\vec{\alpha}(s) = \vec{\alpha}^*(s^*) + \lambda(s^*)\vec{B}^*(s^*) \qquad (2)$$
for some function $\lambda(s^*)$. By taking the derivative of Equation (2) with respect to $s^*$ and using Equations (1), we obtain
$$\vec{T}\frac{ds}{ds^*} = \vec{T}^* + \lambda\tau^*\vec{N}^* + \dot{\lambda}\vec{B}^* \qquad (3)$$
Since $\vec{N}$ and $\vec{B}^*$ are linearly dependent we have $\langle \vec{T}, \vec{B}^* \rangle = 0$. Then, we get
$$\dot{\lambda} = 0.$$
This means that $\lambda$ is a nonzero constant. On the other hand, from the distance function between two points, we have
$$d(\alpha^*(s^*), \alpha(s)) = \|\alpha(s) - \alpha^*(s^*)\|$$
$$= \|\lambda\vec{B}^*\| = |\lambda|\ .$$
Namely, $d(\alpha^*(s^*), \alpha(s)) = \text{constant}$. For the other cases, we obtain the same result.

**Theorem 3.2.** *For a curve $C$ in $E_1^3$, there is a curve $C^*$ such that $\{C, C^*\}$ is a Mannheim pair.*

*Proof.* Since $\vec{N}$ and $\vec{B}^*$ are linearly dependent for all types, Equation (2) can be written as
$$\vec{\alpha}^* = \vec{\alpha} - \lambda\vec{N}\ . \qquad (4)$$
Now, there is a curve $C^*$ for all values of nonzero constant $\lambda$.

**Theorem 3.3.** *Let* $\{C, C^*\}$ *be a Mannheim pair in* $E_1^3$. *Then the relations between the curvatures and torsions of the curves* $C, C^*$ *are given as follows:*

**i)** *If the pair* $\{C, C^*\}$ *is a Mannheim pair of the type 1 or 4, then*

$$\tau^* = \frac{-\kappa}{\lambda \tau}.$$

**ii)** *If the pair* $\{C, C^*\}$ *is a Mannheim pair of the type 2, 3 or 5, then*

$$\tau^* = \frac{\kappa}{\lambda \tau}.$$

*Proof i)* Let the pair $\{C, C^*\}$ be a Mannheim pair of the type 1. By considering the nonzero constant $\lambda$ in Equation (3), we obtain

$$\vec{T} \frac{ds}{ds^*} = \vec{T}^* + \lambda \tau^* \vec{N}^*. \tag{5}$$

Considering Definition 2.1, we have

$$\begin{cases} \vec{T} = \sinh\theta \vec{T}^* + \cosh\theta \vec{N}^* \\ \vec{B} = \cosh\theta \vec{T}^* + \sinh\theta \vec{N}^* \end{cases} \tag{6}$$

where $\theta$ is the angle between the tangent vectors $\vec{T}$ and $\vec{T}^*$ at the corresponding points of the curves $C$ and $C^*$. From Equations (5) and (6), we get

$$\cosh\theta = \lambda \tau^* \frac{ds^*}{ds}, \quad \sinh\theta = \frac{ds^*}{ds} \tag{7}$$

By considering Equation (1), the derivative of Equation (4) with respect to $s^*$ gives us the following

$$\vec{T}^* = (1 - \lambda\kappa)\frac{ds}{ds^*}\vec{T} - \lambda\tau\frac{ds}{ds^*}\vec{B} \tag{8}$$

From Equation (6), we get

$$\begin{cases} \vec{T}^* = -\sinh\theta \vec{T} + \cosh\theta \vec{B} \\ \vec{N}^* = \cosh\theta \vec{T} - \sinh\theta \vec{B} \end{cases} \tag{9}$$

From the Equations (8) and (9), we obtain

$$\cosh\theta = -\lambda\tau\frac{ds}{ds^*}, \quad \sinh\theta = (\lambda\kappa - 1)\frac{ds}{ds^*} \tag{10}$$

Then by the equations (7) and (10), we see that

$$\cosh^2\theta = -\lambda^2 \tau \tau^*, \quad \sinh^2\theta = \lambda\kappa - 1$$

which gives us

$$\tau^* = \frac{-\kappa}{\lambda \tau}.$$

The proof of the statement given in (ii) can be given by a similar way.

**Theorem 3.4.** *Let* $\{C, C^*\}$ *be a Mannheim pair in* $E_1^3$. *The relationship between the curvature and torsion of the curve C is given as follows,*

**i)** *If the pair* $\{C, C^*\}$ *is a Mannheim pair of the type 1, 2 or 5, then we have*

$$\mu\tau + \lambda\kappa = 1$$

**ii)** *If the pair $\{C, C^*\}$ is a Mannheim pair of the type 3, or 4, then the relationship is given by*

$$\mu\tau - \lambda\kappa = 1$$

*where $\lambda$ and $\mu$ are nonzero real numbers.*

*Proof. i)* Assume that the pair $\{C, C^*\}$ is a Mannheim pair of the type 1. Then, from equation (10), we have

$$\frac{-\cosh\theta}{\lambda\tau} = \frac{-\sinh\theta}{1-\lambda\kappa},$$

and so, we get

$$1 - \lambda\kappa = \lambda(\tanh\theta)\tau$$

which gives us

$$\mu\tau + \lambda\kappa = 1$$

where $\lambda$ and $\mu = \lambda\tanh\theta$ are non-zero constants.

The proofs of the statement (ii) can be given by the same way.

**Theorem 3.5.** *Let $\{C, C^*\}$ be a Mannheim pair in $E_1^3$. Then, the relationships between the curvatures and the torsions of the curves $C$ and $C^*$ are given as follows,*

**a)** *If the pair $\{C, C^*\}$ is a Mannheim pair of the type 1, then*

**i)** $\kappa^* = -\dfrac{d\theta}{ds^*}$  **ii)** $\tau^* = \kappa\cosh\theta + \tau\sinh\theta$  **iii)** $\kappa = \tau^*\cosh\theta$  **iv)** $\tau = -\tau^*\sinh\theta$

**b)** *If the pair $\{C, C^*\}$ is a Mannheim pair of the type 2, then*

**i)** $\kappa^* = -\dfrac{d\theta}{ds^*}$  **ii)** $\tau^* = -\kappa\sinh\theta - \tau\cosh\theta$  **iii)** $\kappa = \tau^*\sinh\theta$  **iv)** $\tau = -\tau^*\cosh\theta$

**c)** *If the pair $\{C, C^*\}$ is a Mannheim pair of the type 3, then*

**i)** $\kappa^* = -\dfrac{d\theta}{ds^*}$  **ii)** $\tau^* = -\kappa\sinh\theta + \tau\cosh\theta$  **iii)** $\kappa = \tau^*\sinh\theta$  **iv)** $\tau = \tau^*\cosh\theta$

**d)** *If the pair $\{C, C^*\}$ is a Mannheim pair of the type 4, then*

**i)** $\kappa^* = \dfrac{d\theta}{ds^*}$  **ii)** $\tau^* = \kappa\cosh\theta - \tau\sinh\theta$  **iii)** $\kappa = \tau^*\cosh\theta$  **iv)** $\tau = \tau^*\sinh\theta$

**e)** *If the pair $\{C, C^*\}$ is a Mannheim pair of the type 5, then*

**i)** $\kappa^* = -\dfrac{d\theta}{ds^*}$  **ii)** $\tau^* = \kappa\sin\theta + \tau\cos\theta$  **iii)** $\kappa = \tau^*\sin\theta$  **iv)** $\tau = \tau^*\cos\theta$

*Proof. a)* Let the pair $\{C, C^*\}$ be a Mannheim pair of the type 1 in Minkowski 3-space.

**i)** By taking the derivative of the equation $\langle \vec{T}, \vec{T}^* \rangle = \sinh\theta$ with respect to $s^*$, we have

$$\langle \kappa\vec{N}, \vec{T}^* \rangle + \langle \vec{T}, \kappa^*\vec{N}^* \rangle = \cosh\theta \frac{d\theta}{ds^*}$$

Furthermore, by considering $\vec{N}$ and $\vec{B}^*$ are linearly dependent and using the equations (2) and (9), we have

$$\kappa^* = -\frac{d\theta}{ds^*}.$$

By considering the equations $\langle \vec{N}, \vec{N}^* \rangle = 0$, $\langle \vec{T}, \vec{B}^* \rangle = 0$ and $\langle \vec{B}, \vec{B}^* \rangle = 0$, the proofs of the statements (ii), (iii) and (iv) of (a) in the Theorem 3.5 can be given by the similar way of the proof of the statement (i).

From the statements (iii) and (iv) of the Theorem 3.5, we obtain the following result.

**Proposition 3.1.** *The torsion of the curve $C^*$ is given by*
$$\tau^* = \kappa^2 - \tau^2.$$

The statements (b), (c), (d) and (e) can be proved as given in the proof of the statement (a).

**Theorem 3.6.** *Let $\{C, C^*\}$ be a Mannheim pair in $E_1^3$. For the corresponding points $\alpha(s)$ and $\alpha^*(s^*)$ of the curves $C$, $C^*$ and for the curvature centers $M$ and $M^*$ at this points, that ration*
$$\frac{\|\alpha^*(s^*)M\|}{\|\alpha(s)M\|} : \frac{\|\alpha^*(s^*)M^*\|}{\|\alpha(s)M^*\|}$$
*is not constant.*

*Proof.* Assume that the pair $\{C, C^*\}$ is a Mannheim pair of the type 1. Then, we obtain the following equations,
$$\|\alpha(s)M\| = \frac{1}{\kappa}, \quad \|\alpha^*(s^*)M^*\| = \frac{1}{\kappa^*},$$
$$\|\alpha(s)M^*\| = \sqrt{\left|\lambda^2 - \frac{1}{(\kappa^*)^2}\right|}, \quad \|\alpha^*(s^*)M\| = \frac{1}{\kappa} - \lambda$$

So, we have
$$\frac{\|\alpha^*(s^*)M\|}{\|\alpha(s)M\|} : \frac{\|\alpha^*(s^*)M^*\|}{\|\alpha(s)M^*\|} = (1-\lambda\kappa)\sqrt{\left|\lambda^2(\kappa^*)^2 - 1\right|} \neq \text{constant}.$$

If the pair $\{C, C^*\}$ is a Mannheim pair of the type 2, 3, 4 or 5 we again find that the ration is not constant.

**Proposition 3.2.** *The Mannheim's theorem is invalid for the Mannheim curves in $E_1^3$.*

**Theorem 3.7.** *Let the spherical indicatrix of principal normal vector of the curve $C$ be denoted by $C_2$ with the arclenght parameter $s_2$ and the spherical indicatrix of binormal vector of the curve $C^*$ be denoted by $C_3^*$ with the arclenght parameter $s_3^*$. If $\{C, C^*\}$ is a Mannheim pair in $E_1^3$, then we have the followings,*

i) *If the pair $\{C, C^*\}$ is a Mannheim pair of the type 1, we have*
$$\kappa \frac{ds}{ds_2} = \tau^* \frac{ds^*}{ds_3^*} \cosh\theta, \quad \tau \frac{ds}{ds_2} = -\tau^* \frac{ds^*}{ds_3^*} \sinh\theta.$$

**ii)** *If the pair $\{C, C^*\}$ is a Mannheim pair of the type 2 or 3, we get*

$$\kappa \frac{ds}{ds_2} = -\tau^* \frac{ds^*}{ds_3^*} \sinh \theta, \quad \tau \frac{ds}{ds_2} = -\tau^* \frac{ds^*}{ds_3^*} \cosh \theta.$$

**iii)** *If the pair $\{C, C^*\}$ is a Mannheim pair of the type 4, we have*

$$\kappa \frac{ds}{ds_2} = -\tau^* \frac{ds^*}{ds_3^*} \cosh \theta, \quad \tau \frac{ds}{ds_2} = \tau^* \frac{ds^*}{ds_3^*} \sinh \theta.$$

**iv)** *If the pair $\{C, C^*\}$ is a Mannheim pair of the type 5, we have*

$$\kappa \frac{ds}{ds_2} = \tau^* \frac{ds^*}{ds_3^*} \sin \theta, \quad \tau \frac{ds}{ds_2} = \tau^* \frac{ds^*}{ds_3^*} \cos \theta.$$

*Proof. i)* Suppose that the pair $\{C, C^*\}$ is a Mannheim pair of the type 1. Let $\vec{T}_2$ be the tangent vector of the spherical indicatrix of the principal normal vector of the curve $C$ and $\vec{T}_3^*$ be the tangent vector of the spherical indicatrix of the binormal vector of the curve $C^*$. Since $\vec{N}$ and $\vec{B}^*$ are linearly dependent, the spherical indicatrix of principal normal of the curve $C$ is the same with the spherical indicatrix of binormal of the curve $C^*$. Subsequently, we have

$$\vec{T}_2 = \vec{N}' = (\kappa \vec{T} + \tau \vec{B}) \frac{ds}{ds_2}$$

and

$$\vec{T}_3^* = \vec{B}^{*\prime} = \tau^* \vec{N}^* \frac{ds^*}{ds_3^*}.$$

Since $\vec{N}$ and $\vec{B}^*$ are linearly dependent, we can assume

$$\vec{T}_2 = \vec{T}_3^*$$

Thus, we obtain the following equations

$$\kappa \sinh \theta = -\tau \cosh \theta, \quad \kappa \frac{ds}{ds_2} \cosh \theta + \tau \frac{ds}{ds_2} \sinh \theta = \tau^* \frac{ds^*}{ds_3^*}$$

which gives us

$$\kappa \frac{ds}{ds_2} = \tau^* \frac{ds^*}{ds_3^*} \cosh \theta, \quad \tau \frac{ds}{ds_2} = -\tau^* \frac{ds^*}{ds_3^*} \sinh \theta$$

which are desired equalities.

The proofs of the statements (ii), (iii) and (iv) of the Theorem 3.7, can be given similar way.

**Example 1.** Let consider the spacelike curve $(C^*)$ given by the parametrization

$$\alpha^*(s) = \left( -\frac{1}{2} \sinh s, \frac{1}{2} \cosh s, \frac{\sqrt{5}}{2} s \right).$$

The Frénet vectors of $\alpha^*(s)$ are obtained as follows

$$\vec{T}^* = \left(-\frac{1}{2}\cosh s,\ \frac{1}{2}\sinh s,\ \frac{\sqrt{5}}{2}\right),$$

$$\vec{N}^* = \left(-\sinh s,\ \cosh s,\ 0\right),$$

$$\vec{B}^* = \left(-\frac{\sqrt{5}}{2}\cosh s,\ \frac{\sqrt{5}}{2}\sinh s,\ \frac{1}{2}\right).$$

For $\lambda = 20$, the parametric representation of the Mannheim partner curve $(C)$ of the curve $\alpha^*(s)$ is obtained as

$$\alpha = \left(-\frac{1}{2}\sinh s - 10\sqrt{5}\cosh s,\ \frac{1}{2}\cosh s + 10\sqrt{5}\sinh s,\ \frac{\sqrt{5}}{2}s + 10\right).$$

Then, the pair $\{C, C^*\}$ is a Mannheim pair of the type 3. Figure 1 shows the different appearances of the curves in space in which the curves $\alpha^*$ and $\alpha$ are rendered by blue and red colors, respectively.

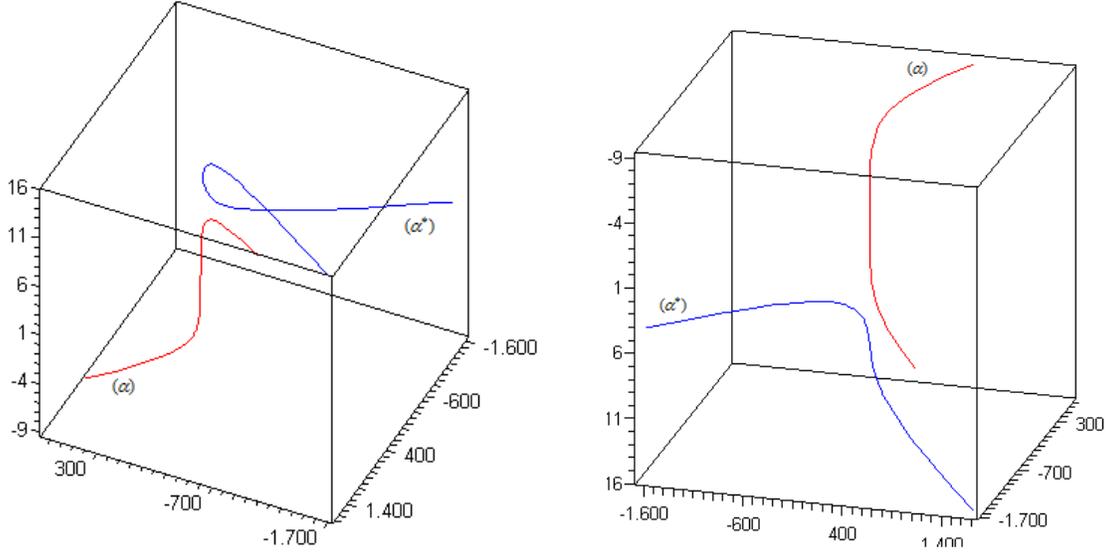

**Fig. 1.** The spacelike curve $\alpha^*$ and its Mannheim partner curve $\alpha$.

**Example 2.** Let now consider the timelike curve $(C^*)$ given by the parametrization
$$\alpha^*(s) = \left(2\sinh s,\ 2\cosh s,\ \sqrt{3}s\right).$$

The Frénet vectors of $\alpha^*(s)$ are obtained as follows

$$\vec{T}^* = \left(2\cosh s,\ 2\sinh s,\ \sqrt{3}\right),$$

$$\vec{N}^* = \left(\sinh s,\ \cosh s,\ 0\right),$$

$$\vec{B}^* = \left(-\sqrt{3}\cosh s,\ -\sqrt{3}\sinh s,\ -2\right).$$

Then for $\lambda = 20$, the Mannheim partner curve $(C)$ of the curve $\alpha^*(s)$ is obtained as

$$\alpha = \left(2\sinh s - 20\sqrt{3}\cosh s,\ 2\cosh s - 20\sqrt{3}\sinh s,\ \sqrt{3}s - 40\right).$$

Then, the pair $\{C, C^*\}$ is a Mannheim pair of the type 1. Figures 2 shows the different appearances of the curves in the space in which the curves $(\alpha^*)(\alpha))$ and $\alpha$ are rendered by blue and red colors, respectively.

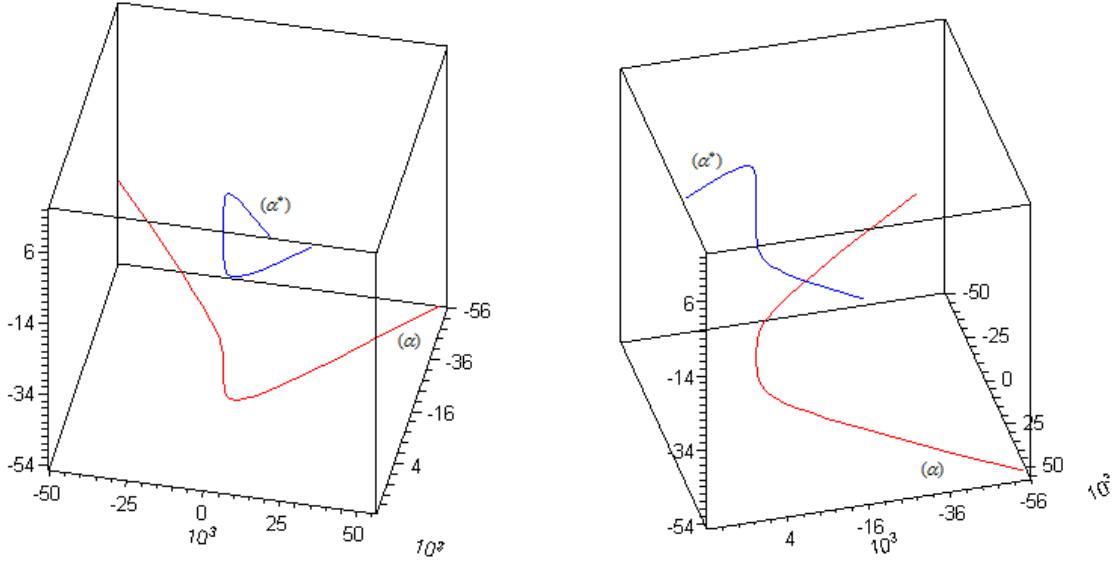

**Fig. 2.** The timelike curve $\alpha^*$ and its Mannheim partner curve $\alpha$.

## 4. CONCLUSIONS

In this paper, we give some characterizations of Mannheim Partner Curves in Minkowski 3-space $E_1^3$. Moreover, we show that the Mannheim theorem is not valid for Mannheim partner curves in $E_1^3$. Also, by considering the spherical indicatrix of some Frénet vectors of the Mannheim curves we give some new characterizations for these curves.

## REFERENCES

1. Burke J. F. Bertrand Curves Associated with a Pair of Curves. *Mathematics Magazine*, 1960, **34**, No. 1, 60-62.
2. Görgülü, E., Ozdamar, E. A generalizations of the Bertrand curves as general inclined curves in $E^n$. *Communications de la Fac. Sci. Uni. Ankara, Series A1*, 1986, **35**, 53-60.
3. Hacisalihoğlu, H. H. *Diferansiyel Geometri.* İnönü Üniversitesi Fen-Edebiyat Fakültesi Yayınları No. 2. 1983.
4. Izumiya, S., Takeuchi, N. Generic properties of helices and Bertrand curves. Journal of Geometry, 2002, **74**, 97-109.
5. Liu, H., Wang, F. Mannheim partner curves in 3-space. *Journal of Geometry*, 2008, **88**, No. 1-2, 120-126.
6. Onder, M., Uğurlu H. H., Kazaz, M. Mannheim offsets of timelike ruled surfaces in Minkowski 3-space. *arXiv:0906.2077v4 [math.DG]*.
7. Onder, M., Uğurlu H. H., Kazaz, M. Mannheim offsets of spacelike ruled surfaces in Minkowski 3-space. *arXiv:0906.4660v3 [math.DG]*.

8. O'Neill, B. *Semi-Riemannian Geometry with Applications to Relativity.* Academic Press, London, 1983.
9. Orbay, K., Kasap, E., Aydemir, I. Mannheim Offsets of Ruled Surfaces. *Mathematical Problems in Engineering*, 2009, Article ID 160917, 9 pages.
10. Orbay, K., Kasap, E. On Mannheim Partner Curves in $E^3$. *International Journal of Physical Sciences,* 2009, 4 (5), 261-264.
11. Oztekin, H.B., Ergüt, M. Null Mannheim curves in the minkowski 3-space $E_1^3$. *Turk J. Math,* 2011, **35**, 107-114.
12. Ravani, B., Ku, T. S. Bertrand Offsets of ruled and developable surfaces. *Comp. Aided Geom. Design,* 1991, **23**, No. 2, 145-152.
13. Struik, D. J. *Lectures on Classical Differential Geometry.* 2nd ed. Addison Wesley, Dover, 1988.
14. Wang, F., Liu, H. Mannheim partner curves in 3-Euclidean space. *Mathematics in Practice and Theory*, 2007, **37**, no. 1, 141-143.
15. Whittemore, J. K. Bertrand curves and helices. *Duke Math. J.* 1940, **6**, No. 1, 235-245.
16. Walrave J. *Curves and surfaces in Minkowski space,* Doctoral thesis, K. U. Leuven, Faculty of Science, Leuven, 1995.